\documentclass[12pt,leqno]{article}
\usepackage[dvips]{graphicx}
\renewcommand{\d}{\displaystyle}
\usepackage{amssymb}
\renewcommand{\d}{\displaystyle}
\begin{document}

\newtheorem{lemma}{Lemma}
\newtheorem{col}[lemma]{Corollary}
\newtheorem{thm}[lemma]{Theorem}
\newtheorem{den}[lemma]{Definition}
\newtheorem{pro}[lemma]{Proposition}
\newcommand{\namelistlabel}[1]{\mbox{#1}\hfil}
\newenvironment{namelist}[1]{%
\begin{list}{}
   {
    \let\makelabel\namelistlabel
    \settowidth{\labelwidth}{#1}
    \setlength{\leftmargin}{1.1\labelwidth}
   }
  }{%
\end{list}}
\renewcommand{\thefootnote}{\fnsymbol{footnote}}
\renewcommand{\theequation}{\arabic{section}.\arabic{equation}}
\def\proof{\noindent{\bf Proof.\quad}}
\def\qed{\hspace*{\fill}\vrule height6pt width4pt depth0pt\medskip}
\def\mod2{\mbox ({\rm mod \hspace{0.2cm} \rm 2})}
\def\frown{\widehat}
\def\ignore#1{\relax}


\normalsize \pagenumbering{arabic} \baselineskip = 24pt

%
\title{Multicolored Isomorphic Spanning Trees in Complete Graphs
\thanks {Research supported in part by NSC 97-2115-M-009-011-MY3.}}

\author{Hung-Lin Fu and Yuan-Hsun Lo \\ \\
{\small \it Department of Applied Mathematics} \\
{\small \it National Chiao Tung University} \\
{\small \it Hsinchu, Taiwan 30050}}

\date{}

\maketitle

\begin{abstract}
In this paper, we first prove that if the edges of $K_{2m}$ are properly colored by $2m-1$ colors in such a way that any two colors induce a 2-factor of which each component is a 4-cycle, then $K_{2m}$ can be decomposed into $m$ isomorphic multicolored spanning trees. Consequently, we show that there exist three disjoint isomorphic multicolored spanning trees in any properly (2$m-$1)-edge-colored $K_{2m}$ for $m\geq 14$.
\end{abstract}

\bigskip
\section{Introduction}

Throughout this paper, all terminologies and notations on graph theory can be referred to the textbook by D. B. West. \cite{W} A {\it spanning tree} $T$ of a graph $G$ is a subgraph of $G$ which is a tree and $V(T)=V(G)$. A {\it $k$-edge-coloring} of a graph $G$ is a mapping from $E(G)$ into a set of colors $\{1,2,\cdots ,k \}$. A $k$-edge-coloring is {\it proper} if incident edges receive distinct colors. Let $\varphi$ be a $k$-edge-coloring of a graph $G$. If $K$ is a subgraph of $G$, for convenience, we use $\varphi|_{{}_K}$ to denote the edge-coloring of $K$ induced by $\varphi$, i.e., $\varphi|_{{}_K}(e)=\varphi(e)$ for each $e\in E(K)$. Note here that in this paper, all edge-colorings are proper.

If $G$ has a proper $k$-edge-coloring, then $G$ is said to be properly $k$-edge-colorable. The {\it chromatic index} of a graph $G$ is the minimum number $k$ such that $G$ is properly $k$-edge-colorable. It's well-known that the chromatic index of $K_{2m}$ is $2m-$1. Let $\varphi$ be a proper (2$m-1$)-edge-coloring of $K_{2m}$ and $C$ be the color set. For each $x\in V(K_{2m})$, define $\varphi_x$ as the mapping from $V(K_{2m})\setminus \{x\}$ to $C$ by $\varphi_x(y)=c$ if $\varphi(xy)=c$. Clearly, $\varphi_x$ is bijective. Let $\varphi^{-1}_x(c)$ be the vertex adjacent to $x$ with the edge colored $c$. For a vertex set $V$ and a color $c$, let $[V]_c=V\cup \{ u|~\varphi(uv)=c, v\in V\}$. For convenience, we use $v\langle c\rangle$ to denote the edge incident to $v$ with color $c$.

A subgraph in an edge-colored graph is said to be {\it multicolored} if no two edges have the same color. Therefore, a question arises naturally: can the edges of a properly (2$m-$1)-edge-colored $K_{2m}$ be partitioned into multicolored subgraphs, such that each has $2m-$1 edges. Here are three conjectures related to this problem.
\vspace{0.7cm}
\\
{\bf Constantine's Conjecture (Weak version) \cite{C}}
\begin{it} For any positive integer m, $m>2$, there exists a proper $(2m-$1$)$-edge-coloring of $K_{2m}$ such that all edges can be partitioned into m isomorphic multicolored spanning trees.
\end{it}
\vspace{0.5cm}
\\
{\bf Brualdi-Hollingsworth Conjecture \cite{BH}}
\begin{it} If $m>2$, then in any proper edge-coloring of $K_{2m}$ with $2m-$1 colors, all edges can be partitioned into m multicolored spanning trees.
\end{it}
\vspace{0.5cm}
\\
{\bf Constantine's Conjecture (Strong version) \cite{C}}
\begin{it} If $m>2$, then in any proper edge-coloring of $K_{2m}$ with $2m-$1 colors, all edges can be partitioned into m isomorphic multicolored spanning trees.
\end{it}

\vspace{0.5cm}

The first conjecture has been proved by Akbari et al \cite{AAFL}. As to the second conjecture, a partial result by Krussel et al \cite{KMV} shows that there are three multicolored spanning trees in $K_{2m}$ for any proper (2$m-$1)-edge-coloring of $K_{2m}$. Essentially, nothing was done so far on the third one. In this paper, we set off the first step by finding three disjoint isomorphic multicolored spanning trees in a proper (2$m-$1)-edge-colored $K_{2m}$ for $m\geq 14$.

It is worth of mention here that the above conjectures will play important roles in applications if they were true. An application of parallelisms of complete designs to population genetics data can be found in \cite{BCML}. Parallelisms are also useful in partitioning consecutive positive integers into sets of equal size with equal power sums \cite{J}. In addition, the discussions of applying colored matchings and design parallelisms to parallel computing appeared in \cite{H}.

\section{The main results}

We start with the notion of a latin square. Let $S$ be an $n$-set. A {\it latin square} of order $n$ based on $S$ is an $n\times n$ array such that each element of $S$ occurs in each row and each column exactly once. For example,
\begin{tabular}{|c|c|}
\hline 0 & 1 \\
\hline 1 & 0 \\
\hline
\end{tabular}
is a latin square of order 2 based on $\{0,1\}=\mathbb{Z}_2$. Since this latin square corresponds to a group table of $\langle \mathbb{Z}_2, + \rangle$, the latin square is also known as a 2-group latin square.

\begin{figure}[h]
    \begin{center}
        \includegraphics[scale=0.95]{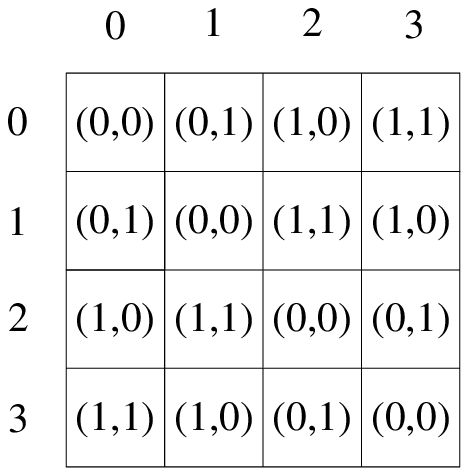}
    \end{center}
    \caption{\label{LxL} $L\times L$.}
\end{figure}

For convenience, a latin square of order $n$ based on $S$ is denoted $L=[~l_{i,j}~]$ where $l_{i,j}\in S$ and $i,j \in \mathbb{Z}_n$. Let $L=[~l_{i,j}~]$ and $M=[~m_{i,j}~]$ be two latin squares of order $l$ and $m$ respectively. Then the direct product of $L$ and $M$ is a latin square of order $l\cdot m$ : $L\times M=[~h_{i,j}~]$ where $h_{x,y}=(~l_{a,b},m_{c,d}~)$ provided that $x=ma+c$ and $y=mb+d$. For instance, let $L$ be the
2-group latin square; then $L\times L$ is a latin square of order 4 based on $\mathbb{Z}_2 \times \mathbb{Z}_2$, as in Figure \ref{LxL}.

A {\it transversal} of a latin square of order $n$ is a set of $n$ entries, one from each column and one from each row, such that these $n$ entries are all distinct. For instance, in $L\times L$, $\{h_{0,0}, h_{1,2}, h_{2,3}, h_{3,1}\}$ is a transversal. $L\times L$ is easily seen to have 4 disjoint transversals. The following shows $L^n=L\times L\times \cdots \times L$ based on ${\mathbb{Z}_2}^n$ has $2^n$ disjoint transversals for each $n\geq 2$.

\begin{pro}
\label{transversal} $L^n$ has $2^n$ disjoint transversals for each $n\geq 2$.
\end{pro}

\proof The proof is by induction on $n$. By Figure 2, $n=2$ is true.

\begin{figure}[h]
    \begin{center}
        \includegraphics[scale=0.8]{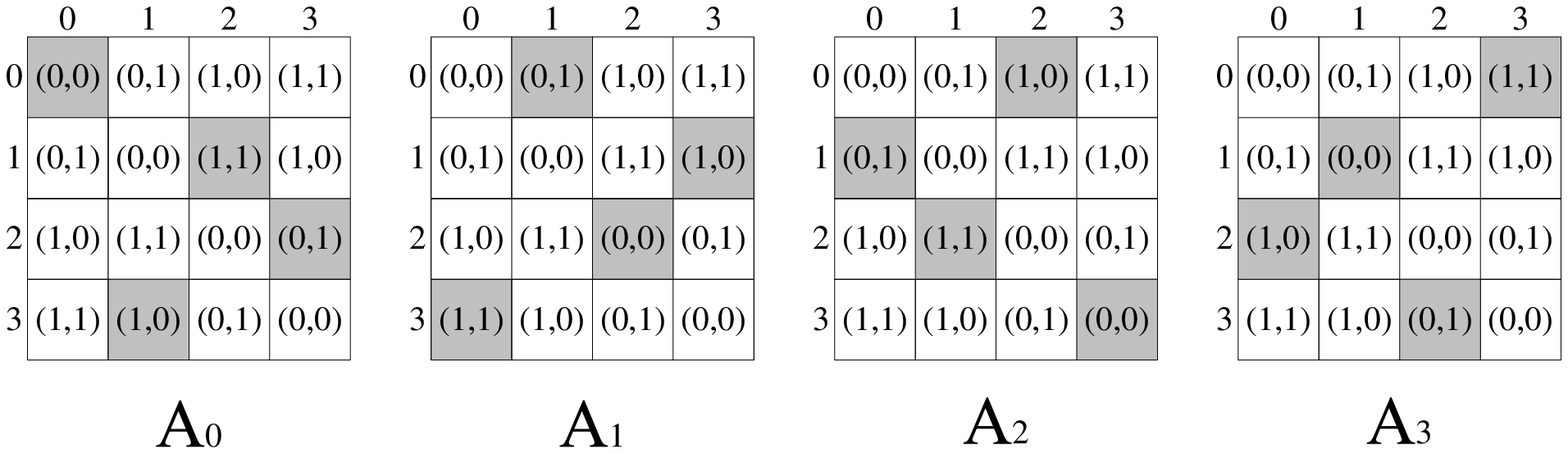}
    \end{center}
    \caption{\label{transversals} 4 transversals in $L^2$.}
\end{figure}

Assume that the assertion is true for each $k\geq 2$. Let $L^k=[{l_{a,b}}^{(k)}]$ and $L^{k+1}=$
\begin{tabular}{|c|c|}
\hline ${L_0}^k$ & ${L_1}^k$ \\
\hline ${L_1}^k$ & ${L_0}^k$ \\
\hline
\end{tabular} .
By definition of direct product, we have ${L_0}^k = [m_{a,b}]$ where $m_{a,b}=(0~,~{l_{a,b}}^{(k)})$ (a ($k+$1)-dim.~vector) and ${L_1}^k = [\overline{m}_{a,b}]$ where $\overline{m}_{a,b}=(1~,~{l_{a,b}}^{(k)})$. We shall use the set of $2^k$ disjoint transversals in $L^k$ to construct $2^{k+1}$ disjoint transversals in $L^{k+1}$.

Let $\{A_i~|~i=0,1,2,\cdots,2^k-1 \}$ be the set of disjoint transversals obtained in $L^k$ by induction hypothesis. Without loss of generality, we may let $A_i$ be the transversal which contains the entry ${l_{0,i}}^{(k)}$, $i=0,1,2,\cdots,2^k-1$. Now, we shall use $A_{2i}$ and $A_{2i+1}$, $i=0,1,2,\cdots,2^{k-1}-1$, to construct four disjoint transversals in $L^{k+1}$. For convenience, we explain the construction by using $A_0$ and $A_1$.

\begin{figure}[h]
    \begin{center}
        \includegraphics[scale=0.8]{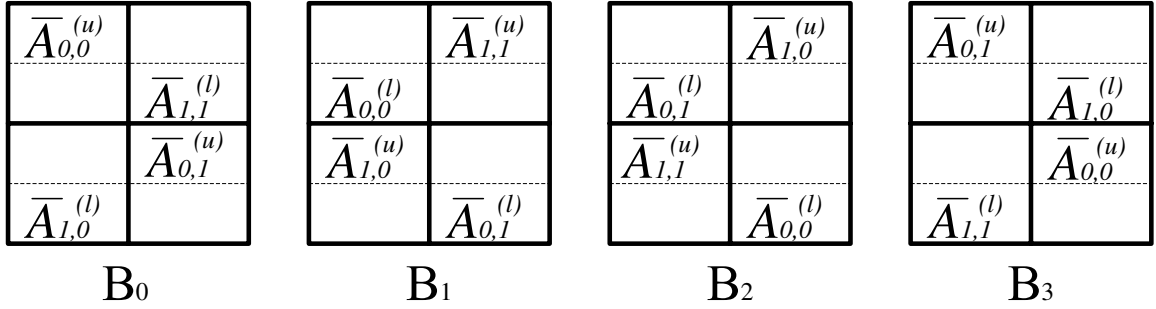}
    \end{center}
    \caption{\label{new_transversals} 4 transversals in $L^{k+1}$ constructed from $A_0$ and $A_1$.}
\end{figure}

Since $A_0$(respectively $A_1$) is a transversal in $L^k$, the corresponding entries in ${L_0}^k$ form a transversal, so are the corresponding entries in ${L_1}^k$. Let the corresponding transversals of $A_0$ in ${L_0}^k$ and ${L_1}^k$ be $\overline{A}_{0,0}$ and $\overline{A}_{1,0}$ respectively. Similarly, let the corresponding transversals of $A_1$ be $\overline{A}_{0,1}$ and $\overline{A}_{1,1}$ respectively. Note that for $0\leq r,s\leq 1$, $\overline{A}_{r,s}$ has $2^k$ entries, one from each row and from each column. Now, for $0\leq r,s\leq 1$, we split $\overline{A}_{r,s}$ into two parts: ${\overline{A}_{r,s}}^{(u)}$ is the set of entries from the first to the $2^{k-1}$-th row of $\overline{A}_{r,s}$, and ${\overline{A}_{r,s}}^{(l)}$ is the set of entries of the other half. By defining $B_0,B_1, B_2$ and $B_3$ as in Figure \ref{new_transversals}, we have four transversals in $L^{k+1}$ as desired.

Since for $i=1,2,\cdots,2^{k-1}-1$, $\overline{A}_{2i}$ and $\overline{A}_{2i+1}$ can also be used to construct four transversals in $L^{k+1}$, we have a set of $2^{k+1}$ transversals in $L^{k+1}$. By the reason that $A_0, A_1,\cdots, A_{2^k-1}$ are disjoint transversals, we conclude the proof. \qed

\begin{lemma} \label{unique} If $\mu$ is a proper $(2m-$1$)$-edge-coloring of $K_{2m}$, $m\geq 2$, such that any two colors induce a 2-factor with each component a 4-cycle, then {\rm (a)} $2m=2^n$ for some $n\geq 2$ and {\rm (b)} $K_{2m}$ contains a clique $K$ of order $2^k$, $1\leq k \leq n-1$ such that $\{\mu(e)~|~e\in E(K)\}$ is a $(2^k-$1$)$-set, i.e., $\mu|_{{}_K}$ is a $(2^k-$1$)$-edge-coloring of $K$.
\end{lemma}

\proof First, we claim that (b) is true. The proof is by induction on $n$. Clearly, it is true when $n=2$. By hypothesis, let $H$ be a clique of order $2^h$, $h<k$, and $\mu|_{{}_H}$ is a $(2^h-$1$)$-edge-coloring of $H$. Without loss of generality, let $V(H)=\{x_1,x_2,\cdots,x_{2^h}\}$ and the colors used in $H$ be $\{c_1,c_2,\cdots,c_{2^h-1}\}$. Since $\mu$ is a proper (2$m-$1)-edge-coloring of $K_{2m}$, each color occurs around each vertex. Let $c_{2^h}$ be a color not used in $H$. Then, we have a set $H'$, $H'\cap H=\phi$, $H'=\{y_1,y_2,\cdots,y_{2^h}\}$ such that $\mu(x_iy_i)=c_{2^h}$ for $i=1,2,\cdots,2^h$. Now, by the reason that any two colors induce a $C_4$-factor, we conclude that $\mu|_{{}_{H'}}$ is also a $(2^h-$1$)$-edge-coloring of $H'$, moreover, $\mu(x_ix_j)= \mu(y_iy_j)$ for $1\leq i \neq j \leq 2^h$. Therefore, the complete bipartite graph $K_{2^h,2^h}=(H,H')$ has a $2^h$-edge-coloring following by the same reason. This implies that $\mu|_{{}_{H\cup H'}}$ is a $(2^{h+1}-$1$)$-edge-coloring of the clique induced by $H\cup H'$. So, we have the proof of (b).

Suppose $2m=2^r\cdot p$ where $p$ is an odd integer and $p\neq 1$. Using the above argument, we can find the largest clique $G$ of order $2^s$ which uses $2^s-1$ colors. Then we partition the vertices of $K_{2m}$ into two sets $X$ and $Y$ where $X=V(G)$, and let $|Y|=q$. Here, we notice that $q<2^s$. Consider these $2^s-1$ colors used in coloring the edges of $G$, there are total $(2^s-1)(2^{r-1}\cdot p)$ edges which use these colors. But, we have used these colors in $G$. Hence, there remains $\d \frac12(2^s-1)(2m-2^s)$ edges to be colored by using these colors. Since the edges between X and Y can't be colored with any of these colors, they have to be in Y. But, since $q<2^s$ and $2m-2^s=q$, $\d \frac12(2^s-1)(2m-2^s) > {q \choose 2}$, a contradiction. This implies that $p=1$, and we have the proof of (a). \qed

\begin{lemma} {\rm\cite{BH}} \label{partitionK8} Let $\mu$ be a proper 7-edge-coloring of $K_8$ such that for any two colors form a $C_4$-factor. Then the edges of $K_8$ can be partitioned into 4 isomorphic multicolored spanning trees.
\end{lemma}

\begin{thm} \label{partition} If $\mu$ is a proper $(2m-$1$)$-edge-coloring of $K_{2m}$, $m>2$, such that any two colors form an $C_4$-factor, then the edges of $K_{2m}$ can be partitioned into m isomorphic multicolored spanning trees.
\end{thm}

\proof By Lemma \ref{unique}, $2m=2^n$ for some $n>2$. We prove the theorem by induction on $n$. By Lemma \ref{partitionK8}, $n=3$ is true.

Assume that the assertion is true for each $k\geq 3$ and consider $K_{2^{k+1}}$.

From the process of the proof of Lemma \ref{unique}, there must exist two disjoint cliques of order $2^k$ with $2^k-$1 colors in $K_{2^{k+1}}$. Let $V(K_{2^{k+1}})=A\cup B$ where $A, B$ are the vertex sets of the two cliques. Consider the colors of the edges between $A$ and $B$. Let $A=\{a_0, a_1, \ldots, a_{2^k-1}\}$, $B=\{b_0, b_1, \ldots, b_{2^k-1}\}$ and $M=[m_{i,j}]$ where $m_{i,j}= \mu(a_ib_j)$. It's clear that $M$ is a latin square; furthermore, $M\cong L^k$. By Proposition \ref{transversal}, $M$ has $2^k$ disjoint transversals. This implies that there are $2^k$ perfect matchings in the complete bipartite graph induced by $A\cup B$. Note that the two cliques induced by $A$ and $B$ respectively have $2^{k-1}$ multicolored isomorphic spanning trees of order $2^k$, respectively. Thus, by assigning a perfect matching to each spanning tree, we obtain $2^k$ spanning trees of order $2^{k+1}$. Moreover, these spanning trees are isomorphic and multicolored. \qed

For the presentation of the proof of our main theorem, we define the following notations. In a properly ($2m-$1)-edge-colored $K_{2m}$, a {\it u-star} $S_u$ is a spanning tree consisting of all edges incident to $u$, where $u\in V(K_{2m})$. Suppose $T$ is a multicolored spanning tree of $K_{2m}$ with two leaves $x_1$ and $x_2$. Let the edges incident to $v_1$ and $v_2$ be $e_1$ and $e_2$ respectively, and $\varphi(e_1)=c_1$, $\varphi(e_2)=c_2$. Then let $T[x_1,x_1;c_1,c_2]$ be the tree obtained from $T$ by removing the edges $e_1, e_2$ and adding the edges $x_1\langle c_2\rangle, x_2\langle c_1\rangle$.

At first, we show the existence of two disjoint isomorphic multicolored spanning trees.

\begin{lemma} \label{two trees} Let $\varphi$ be an arbitrary proper $(2m-$1$)$-edge-coloring of $K_{2m}$. Then there exist two disjoint isomorphic multicolored spanning trees in $K_{2m}$ for $m\geq 3$.
\end{lemma}

\proof Let $V(K_{2m})=\{x_i| ~i=1,2,\ldots,2m \}$. We split the proof into two cases.

\begin{description}

\item[Case 1.] There exists a 4-cycle $(x_1, x_2, x_3, x_4)$ such that $\varphi(x_1x_2)=b$, $\varphi(x_3x_4)=c$, and $\varphi(x_1x_4)$ $=\varphi(x_2x_3)=a$. Let $T_1=S_{x_1}[x_2,x_4;b,a]$ and $T_2=S_{x_2}[x_1,x_3;b,a]$, see Figure \ref{2_trees}. Clearly, they are the desired spanning trees.

\begin{figure}[h]
    \begin{center}
        \includegraphics[scale=0.8]{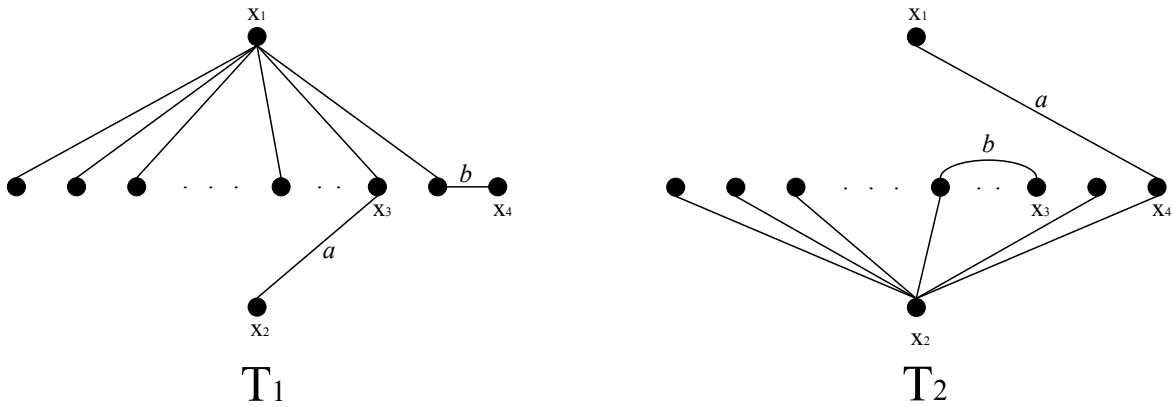}
    \end{center}
    \caption{\label{2_trees} Two isomorphic spanning trees of Case 1.}
\end{figure}

\item[Case 2.] If any two colors of this edge-coloring induce a $C_4$-decomposition of $K_{2m}$, then we have the proof by Theorem \ref{partition}. \qed
\end{description}

Now, we are ready for the main result.

\begin{thm} \label{three trees} Let $\varphi$ be an arbitrary proper $(2m-$1$)$-edge-coloring of $K_{2m}$. Then there exist three disjoint isomorphic multicolored spanning trees in $K_{2m}$ for $m\geq 14$. \end{thm}

\proof From the proof of Lemma \ref{two trees}, we only need to consider the case: there exist two colors which do not induce a $4$-cycle factor. Let $T_1$ and $T_2$ be the isomorphic multicolored spanning trees obtained in Lemma \ref{two trees}. Clearly, $K_{2m}-T_1-T_2$ is disconnected ($\{x_1,x_2\}$ induces a component in this graph). Let $\varphi^{-1}_{x_3}(b)=y_1$, $\varphi^{-1}_{x_4}(b)=y_2$ and $U=V(K_{2m})-\{x_1,x_2,x_3,x_4,y_1,y_2\}$. Since $m\geq 14$, we can choose a vertex $u\in U$ such that the two colors $\varphi(ux_1)$ and $\varphi(ux_2)$ are different from those colors on the edges of the graph induced by the vertex set $\{x_1, x_2, x_3, x_4\}$. Without loss of generality, let $\varphi(ux_1)=1$ and $\varphi(ux_2)=2$. Moreover, let $v_1\in U~\backslash ~\{u\}$ and  $\varphi(x_1v_1)=3$ such that $\varphi^{-1}_{v_1}(b)\neq \varphi^{-1}_{x_4}(1)$ and the two vertices $\varphi^{-1}_u(3)$ and $\varphi^{-1}_{v_1}(1)$ are elements in $U~\backslash ~\{u\}$. Now, pick $v_2\in U~\backslash ~\{u,v_1,\varphi^{-1}_{v_1}(b) \}$ and let $\varphi(x_2v_2)=4$ such that $\varphi^{-1}_{v_2}(b)\neq \varphi^{-1}_{x_3}(2)$ and the two vertices $\varphi^{-1}_u(4)$ and $\varphi^{-1}_{v_2}(2)$ are elements in set $U~\backslash ~\{u\}$. Note that we can always pick $v_1$ and $v_2$ consecutively since $m\geq 14$.

\begin{figure}[h]
    \begin{center}
        \includegraphics[scale=0.8]{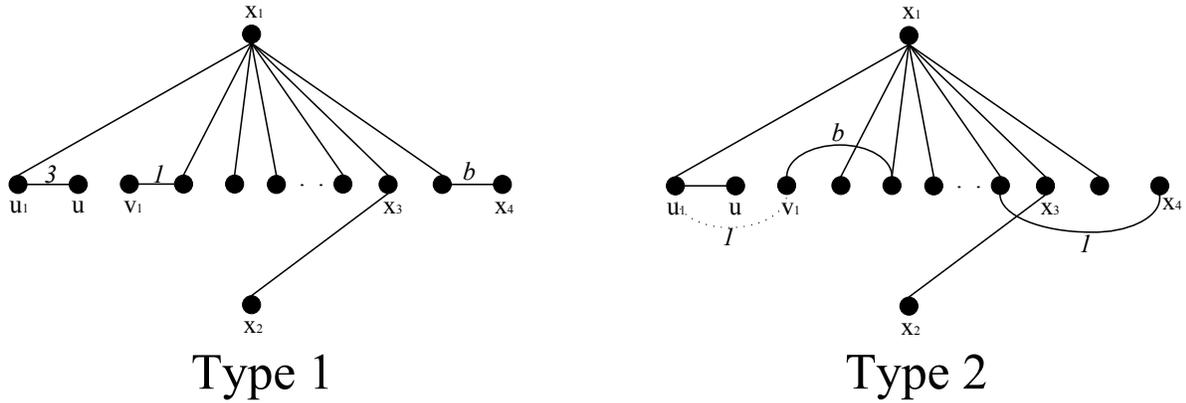}
    \end{center}
    \caption{\label{T1} Two types of $T_1'$.}
\end{figure}

Let $T_1'=T_1[u,v_1;1,3]$ and $T_2'=T_2[u,v_2;2,4]$. Assume that $\varphi^{-1}_u(3)=u_1$ and $\varphi^{-1}_u(4)=u_2$. If $u_1=\varphi^{-1}_{v_1}(1)$, then adjust $T_1'$ to $T_1'[v_1,x_4;1,b]$. Similarly, if $u_2=\varphi^{-1}_{v_2}(2)$, then adjust $T_2'$ to $T_2'[v_2,x_3;2,b]$. Then $T_1'$ and $T_2'$ both have two types. In either case, they are disjoint and isomorphic. Figure \ref{T1} shows the types of $T_1'$.

\begin{figure}[h]
    \begin{center}
        \includegraphics[scale=1]{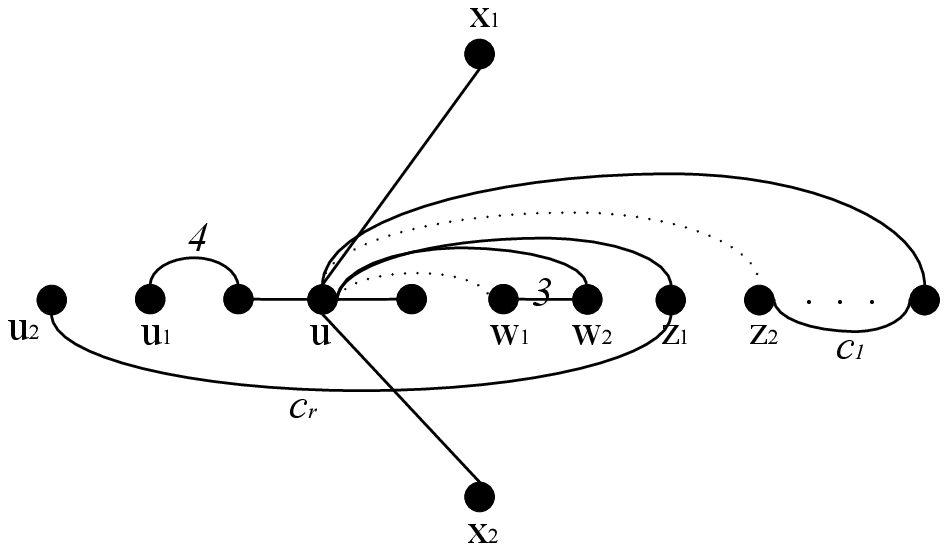}
    \end{center}
    \caption{\label{T3} $T_3'$.}
\end{figure}

Now, we are ready to construct the third tree. Let $T_3$ be the graph $S_u[u_1,u_2;3,4]$. Then choose one edge $w_1w_2$ with color $3$ in the graph induced by $V(K_{2m})\backslash \{x_1,x_2,u,u_2\}$ and assume $\varphi(uw_1)=c_1$, $\varphi(uw_2)=c_2$. Let $W=\{x_1,x_2,u_1, \varphi^{-1}_{u_1}(4),w_1,w_2\}$. Since $m\geq 14$, there exists one color, $c_r$, such that $\varphi^{-1}_{u_2}(c_r)\notin W$ and $\varphi^{-1}_u(c_r)\notin [W]_{c_1} \cup [W]_{c_2}$. Let $\varphi^{-1}_{u_2}(c_r)=z_1$ and $\varphi^{-1}_u(c_r)=z_2$. Since $\varphi(z_1z_2)$ may be $c_1$ or $c_2$, we can assume $\varphi(z_1z_2)\neq c_1$. Finally, let $T_3'$ be obtained from $T_3$ by removing the edges $u_2\langle 3\rangle,u\langle c_p\rangle, u\langle c_r\rangle$ and then adding the edges $u_2\langle c_r\rangle,w_1\langle 3\rangle,z_2\langle c_p\rangle$. Thus, the third spanning tree is constructed, see Figure \ref{T3}. Since all spanning trees contain exactly four vertices which are of distance 2 from vertices $x_1,x_2$ and $u$ respectively, they are isomorphic. This concludes the proof. \qed

\noindent{\bf \Large Acknowledgements}

The authors would like to express their gratitude to the referee for his careful reading and his important comments that significantly improve the presentation of this paper.

\rm

\end{document}